\documentclass[12pt, letterpaper, reqno]{amsart}
\usepackage[hyperindex,breaklinks]{hyperref}
\usepackage{amsmath}
\usepackage{amsfonts}
\usepackage{mathrsfs}
\usepackage{graphicx}
\usepackage{color}
\usepackage{slashed}
\usepackage{braket}
\usepackage{extarrows}
\usepackage{amsmath, bm}
\usepackage{amssymb}
\usepackage[margin=1in]{geometry}
\usepackage{bbm}
\usepackage{hyperref}

\definecolor{yellow}{rgb}{1,.7,0}

\definecolor{pkured}{rgb}{0.55,0,0}

\newcommand{\be}{\begin{equation*}}
	\newcommand{\ee}{\end{equation*}}
\newcommand{\beq}{\begin{equation}}
	\newcommand{\eeq}{\end{equation}}

\numberwithin{equation}{section}

\newtheorem{cor}{Corollary}[section]

\newtheorem{lem}[cor]{Lemma}

\newtheorem{thm}[cor]{Theorem}

\newtheorem{ex}[cor]{Example}
\theoremstyle{remark}
\newtheorem{rmk}[cor]{Remark}

\numberwithin{figure}{section}
\usepackage{ytableau}

\usepackage{tikz}
\newcounter{x}
\newcounter{y}
\newcounter{z}

\author{Chenglang Yang}
\email{yangcl@pku.edu.cn}
\address{Institute for Math and AI, Wuhan University, Wuhan 430072, China}
\address{Hua Loo-Keng Center for Mathematical Sciences \&
	Academy of Mathematics and Systems Science,
	Chinese Academy of Sciences,
	Beijing 100190, China}

\title[{Hook-Lengths, Symplectic/Orthogonal Content, Amdeberhan's Conjectures}]
{Hook-Lengths, Symplectic/Orthogonal Contents and Amdeberhan's Conjectures}

\begin{document}
\maketitle

\begin{abstract}
	The symplectic/orthogonal contents of partitions are related to the dimensions of irreducible representations of symplectic/orthogonal groups.
	In 2012, motivated by Nekrasov--Okounkov's hook-length formula and Stanley's hook-content formula,
	Amdeberhan proposed several conjectures about infinite product formulas for certain generating functions of hook-lengths and symplectic/orthogonal contents.
	Some special cases of his conjectures were recently proved by Amdeberhan, Andrews and Ballantine.
	In this paper,
	we prove the general cases of Amdeberhan's conjectures.
\end{abstract}

\setcounter{section}{0}
\setcounter{tocdepth}{2}


\renewcommand{\thefootnote}{\fnsymbol{footnote}} 
\renewcommand{\thefootnote}{\arabic{footnote}}

\section{Introduction}

The study of integer partitions was initiated by Euler who showed the generating function of the numbers of partitions as an infinite product formula.
Later,
mathematicians intensively studied various aspects of integer partitions,
including their relations to combinatorics, representation theory, and mathematical physics.
For example,
the irreducible representations of symmetric groups,
as well as the polynomial representations of general linear groups,
can be indexed by partitions.
The dimension of representations can be computed by using the hook-lengths of corresponding partitions.
In 2006,
Nekrasov and Okounkov \cite{NO06} provided a formula when they studied the Seiberg--Witten theory and random partitions. This formula relates a sum over products of hook-lengths of partitions to the powers of Euler's infinite product formula.
Various new proofs and generalizations of this Nekrasov--Okounkov formula have been explored,
as discussed in \cite{H08,H10,INRS,S10,W06,Y23b} and references therein.
Inspired by Nekrasov--Okounkov's hook-length formula,
and Stanley's analogous hook-content formula \cite{S10},
Amdeberhan proposed several conjectures about infinite product formulas for certain generating functions of hook-lengths and symplectic/orthogonal contents \cite{Amd12}.
These conjectures may have interests in the representation theory of symplectic/orthogonal groups
(see \cite{SK79,Lit,Sun90b})
and mathematical physics, especially in integrable systems (see, for examples, \cite{JLW,SWC}).
Recently,
some special cases of Amdeberhan's conjectures were proved by Amdeberhan, Andrews and Ballantine \cite{AAB}.
In this paper,
we will prove the general cases of Amdeberhan's conjectures.
The key point is the use of vertex operators.

A partition of a nonnegative integer $n$ is a sequence of weakly decreasing positive integers $\lambda=(\lambda_1,...,\lambda_l)\vdash n$ satisfying $n=|\lambda|=\sum_{i=1}^l \lambda_i$,
where $l$ is called its length.
In general,
one can also write $\lambda_i=0$ for $i>l$.
The conjugation of $\lambda$ is a new partition $\lambda^t$ defined by
$\lambda^t_i=\#\{j|\lambda_j\geq i\}$ for $i=1,...,\lambda_1$.
For each pair of positive integers $(i,j)\in\lambda$,
i.e., $1\leq j\leq \lambda_i$,
one can define the following combination numbers related to the given partition $\lambda$.
The content is $c(i,j):=j-i$ and the hook-length is $h(i,j):=\lambda_i+\lambda^t_j-i-j+1$.
We recommend the books \cite{And98,Mac} to readers for more interesting properties of partitions.
The Amdeberhan's conjectures concern the so called symplectic contents and orthogonal contents of partitions,
which are defined as (see \cite{Sun90b})
\begin{align}
	c_{sp}(i,j)=\begin{cases}
		\lambda_i+\lambda_j-i-j+2,\quad\quad& \text{\ if\ } i>j,\\
		i+j-\lambda^t_i-\lambda^t_j,\quad\quad& \text{\ if\ } i\leq j,
	\end{cases} \label{eqn:def csp}\\
	c_{O}(i,j)=\begin{cases}
		\lambda_i+\lambda_j-i-j,\quad\quad& \text{\ if\ } i\geq j,\\
		i+j-\lambda^t_i-\lambda^t_j-2,\quad\quad& \text{\ if\ } i<j.
	\end{cases} \label{eqn:def co}
\end{align}
These symplectic/orthogonal contents and the hook-lengths are related to the representations of symplectic/orthogonal groups.
Therefore,
it is also interesting to explain the results of this paper in the context of representation theory.
Interested readers may find \cite{And98,S99,Sun90b} useful.

Stanley \cite{S10},
inspired by the Nekrasov--Okounkov formula \cite{NO06},
proposed the hook-content formula,
which relates a sum over products of hook-lengths and contents to a product formula (see Theorem 2.2 in \cite{S10}).
The hook-contents are related to the representation theory of general linear groups.
Thus it is natural to consider its symplectic and orthogonal counterparts,
which should involve symplectic contents and orthogonal contents.
These are exactly the Amdeberhan's conjectures in \cite{Amd12}.

Our first result is stated as follows
\begin{thm}[=Conjecture 6.2 (a) in \cite{Amd12}]
	\label{thm:main csp}
	For indeterminates $t, q$,
	we have
	\begin{align}\label{eqn:main csp}
		\sum_{\lambda\in\mathcal{P}}q^{|\lambda|}
		\prod_{u\in\lambda}
		\frac{t+c_{sp}(u)}{h(u)}
		=\prod_{n=1}^\infty
		\frac{(1-q^{8n})^{\binom{t+1}{2}}}
		{(1-q^{8n-2})^{\binom{t+1}{2}-1}}
		\Big(\frac{1-q^{4n-1}}{1-q^{4n-3}}\Big)^t
		\Big(\frac{1-q^{8n-4}}{1-q^{8n-6}}\Big)^{\binom{t}{2}-1},
	\end{align}
	where $\mathcal{P}$ is the set of partitions.
\end{thm}

It is known that after taking conjugations of partitions,
the symplectic contents and the orthogonal contents are dual to each other.
Based on this combinatorial fact,
we prove
\begin{cor}[=Conjecture 6.2 (b) in \cite{Amd12}]
	\label{cor:main co}
	We have
	\begin{align}\label{eqn:main co}
		\sum_{\lambda\in\mathcal{P}}q^{|\lambda|}
		\prod_{u\in\lambda}
		\frac{t+c_{O}(u)}{h(u)}
		=\prod_{n=1}^\infty
		\frac{(1-q^{8n})^{\binom{t}{2}}}
		{(1-q^{8n-6})^{\binom{t}{2}-1}}
		\Big(\frac{1-q^{4n-1}}{1-q^{4n-3}}\Big)^t
		\Big(\frac{1-q^{8n-4}}{1-q^{8n-2}}\Big)^{\binom{t+1}{2}-1}.
	\end{align}
\end{cor}

Some special cases of the above two results are independently interesting.
For example, Conjecture 6.2 (c) in \cite{Amd12} is the special case $t=0$ of the above two formulas and will be used in Corollary \ref{cor:main csp^2=csp}. 
Next,
after establishing a Cauchy-identity type formula for symplectic Schur functions,
we prove the following
\begin{thm}
	\label{thm:main csp^2 cO^2}
	For indeterminates $t_1, t_2, q$,
	we have
	\begin{align}
		\sum_{\lambda\in\mathcal{P}}q^{|\lambda|}
		\prod_{u\in\lambda}
		\frac{\big(t_1+c_{sp}(u)\big)\big(t_2+c_{sp}(u)\big)}{h^2(u)}
		=\prod_{n=1}^\infty
		\frac{(1-q^{4n-2})^{\binom{t_1}{2}+\binom{t_2}{2}-1}
			(1-q^{4n})^{\binom{t_1+1}{2}+\binom{t_2+1}{2}}}
		{(1-q^{4n-3})^{t_1t_2} (1-q^{4n-1})^{t_1t_2}},
		\label{eqn:main spsp}\\
		\sum_{\lambda\in\mathcal{P}}q^{|\lambda|}
		\prod_{u\in\lambda}
		\frac{\big(t_1+c_{O}(u)\big)\big(t_2+c_{O}(u)\big)}{h^2(u)}=\prod_{n=1}^\infty
		\frac{(1-q^{4n-2})^{\binom{t_1+1}{2}+\binom{t_2+1}{2}-1}
			(1-q^{4n})^{\binom{t_1}{2}+\binom{t_2}{2}}}
		{(1-q^{4n-3})^{t_1t_2} (1-q^{4n-1})^{t_1t_2}}.
		\label{eqn:main coco}
	\end{align}
\end{thm}
Conjecture 6.3 (a) in \cite{Amd12} follows from setting $t_1=-t_2$ and replacing $t_2^2$ by $-t$ in Theorem \ref{thm:main csp^2 cO^2}.
See Example \ref{exa:conj6.3a} for more details.
As a consequence of Theorem \ref{thm:main csp} and Theorem \ref{thm:main csp^2 cO^2}, we have
\begin{cor}[=Conjecture 6.3 (c) in \cite{Amd12}]
	\label{cor:main csp^2=csp}
	For any nonnegative integer $n$,
	we have
	\begin{align}
		(-1)^{\binom{n}{2}}
		\cdot\sum_{\lambda\vdash n}
		\prod_{u\in\lambda}
		\frac{c^2_{sp}(u)}{h^2(u)}
		=\sum_{\lambda\vdash n}
		\prod_{u\in\lambda}
		\frac{c_{sp}(u)}{h(u)}.
	\end{align}
\end{cor}
The cases $t=0$ of Theorem \ref{thm:main csp} and Corollary \ref{cor:main co},
and the cases $t_1=-t_2=0$ and $t_1=-t_2=1$ of Theorem \ref{thm:main csp^2 cO^2} were recently proved by Amdeberhan, Andrews and Ballantine in \cite{AAB}.
Using the method developed in this paper,
and Sundaram's Cauchy identity for symplectic Schur functions (see Theorem 15.1 in \cite{Sun86}),
we can also provide the following new formula,
which involves symplectic contents and ordinary contents at the same time.
\begin{thm}\label{thm:main csp*c}
	We have
	\begin{align}\label{eqn:main csp*c}
		\sum_{\lambda\in\mathcal{P}}q^{|\lambda|}
		\prod_{u\in\lambda}
		\frac{\big(t_1+c_{sp}(u)\big)\big(t_2+c(u)\big)}{h^2(u)}
		=\frac{(1-q^2)^{\binom{t_2}{2}}}{(1-q)^{t_1t_2}}.
	\end{align}
\end{thm}

The rest of this paper is organized as follows.
In Section \ref{sec:pre},
we review the concepts of Schur functions, symplectic Schur functions and vertex operators.
The latter serves as the main tool in this paper.
In Section \ref{sec: proof single},
by studying a partition function consisting of symplectic Schur functions,
we prove Theorem \ref{thm:main csp} and Corollary \ref{cor:main co}.
In Section \ref{sec: proof double},
after establishing a Cauchy-identity type formula for symplectic Schur functions,
we prove Theorems \ref{thm:main csp^2 cO^2}, \ref{thm:main csp*c} and Corollary \ref{cor:main csp^2=csp}.

\section{Preliminaries}
\label{sec:pre}
In this section,
we review the concepts of Schur functions, symplectic Schur functions and the vertex operators.
We recommend the Macdonald's book \cite{Mac} and \cite{DJM}.

\subsection{Schur functions and symplectic Schur functions}
The theory of symmetric functions plays an important role in the representation theory.
In general,
one starts from a set of independent variables $(x_1,...,x_n)$ and denote by $\Lambda_n(x_1,...,x_n)$ the set of symmetric polynomial functions in these variables.
It is well known that the following power sum coordinates
\begin{align}\label{eqn:def p}
	p_k(x_1,...,x_n)=\sum_{i=1}^n x_i^k
\end{align}
form an algebraic basis for this ring.
Thus,
one can write $\Lambda_n(x_1,...,x_n)=\mathbb{C}[\bm p]$ if there is no ambiguity,
where $(\bm p)=(p_1,p_2,...)$.

It is more convenient to extend the notion of symmetric functions to functions in infinitely many variables $(\bm x)=(x_1,x_2,...)$.
In terms of certain projective limit,
one can obtain the ring $\Lambda(\bm x)$ consisting of symmetric functions in infinitely many variables $(\bm x)=(x_1,x_2,...)$
(for more details, see Chapter I.2 in \cite{Mac}).
Here we are focusing on some concrete symmetric functions.
The $n$-th power sum coordinates $p_n(\bm x)$ are still defined by equation \eqref{eqn:def p}, but they now involve infinitely many variables $(\bm x)$.
General symmetric functions could also be represented as polynomials in terms of these power sum coordinates.
For example,
the famous Schur functions can be defined in this manner by using the generating function of elementary Schur functions and the Jacobi--Trudi formula.
In this paper,
we formally define them as follows.
The Schur function labeled by a partition $\lambda=(\lambda_1,...,\lambda_l)$ is
\begin{align*}
	s_\lambda
	=s_\lambda(\bm x)
	=\lim_{n\rightarrow\infty}
	\frac{\det(x_j^{\lambda_i+n-i})_{1\leq i,j\leq n}}
	{\det(x_j^{n-i})_{1\leq i,j\leq n}}.
\end{align*}
A kind of generalization of Schur functions is the skew Schur functions.
They are labeled by two partitions.
More concretely,
for two given partitions $\lambda, \mu$,
the skew Schur function labeled by $(\lambda, \mu)$ is defined by
\begin{align}\label{eqn:def skew schur}
	s_{\lambda/\mu}=\sum_{\nu\in\mathcal{P}} c^{\lambda}_{\mu \nu} s_{\nu},
\end{align}
where $\{c^{\lambda}_{\mu \nu}\}$ are the Littlewood-Richardson coefficients defined by $s_{\mu} s_{\nu}=\sum_{\lambda} c^{\lambda}_{\mu \nu} s_{\lambda}$.

To introduce the symplectic Schur functions,
one should consider symmetric functions involving $2n$ variables $(x_1,x_1^{-1},...,x_n,x_n^{-1})$.
Then,
for any partition $\lambda$ satisfying $l(\lambda)\leq n$,
the symplectic Schur function labeled by $\lambda$ and $n$ is defined by
\begin{align*}
	sp_{\lambda,2n}(x_1,x_1^{-1},...,x_n,x_n^{-1})
	=
	\frac{\det(x_j^{\lambda_i+n-i}-x_j^{-\lambda_i-n+i-1})_{1\leq i,j\leq n}}
	{\det(x_j^{n-i+1}-x_j^{-n+i-1})_{1\leq i,j\leq n}}.
\end{align*}
It is also known that the symplectic Schur function admits the following Jacobi-Trudi formula
(see \cite{KT87,Weyl}, and see \cite{FK97} for a bijective proof)
\begin{align*}
	sp_{\lambda,2n}(x_1^\pm,...,x_n^\pm)
	=\frac{1}{2}\det\big(h_{\lambda_{i}-i+j}(x_1^\pm,...,x_n^\pm)
	+h_{\lambda_i-i-j+2}(x_1^\pm,...,x_n^\pm)\big)_{i,j=1}^n,
\end{align*}
where $(x_1^\pm,...,x_n^\pm)=(x_1,x_1^{-1},...,x_n,x_n^{-1})$,
and $h_k(\cdot)$ is the $k$-th complete homogeneous symmetric function in $2n$ variables.
Thus, the symplectic Schur function $sp_{\lambda,2n}$ could also be regarded as a symmetric function in the variables $(x_1,x_1^{-1},...,x_n,x_n^{-1})$.
In dealing with infinitely many variables,
we adopt the following notation
\begin{align*}
	(\tilde{\bm x})=(\bm x, \bm x^{-1})
	=(x_1,x_1^{-1},...,x_n,x_n^{-1},...).
\end{align*}
Then the symplectic Schur function $sp_\lambda(\tilde{\bm x})$ is just formally defined by the limit $n\rightarrow\infty$ of $sp_{\lambda,2n}$.
In this case,
$sp_\lambda(\tilde{\bm x})$ could also be represented by a polynomial in those power sum coordinates ${p_k(\tilde{\bm x})}$ and we write it as $sp_\lambda(\bm p)$,
where $(\bm p)=(p_1,p_2,...)$, if there is no ambiguity.

\begin{rmk}\label{rmk:everything as p_k}
	It is worth mentioning that,
	when using the notation $sp_{\lambda,2n}$,
	there is a hidden condition $l(\lambda)\leq n$.
	Thus,
	when considering properties involving symplectic Schur functions labeled by all partitions,
	it is better to use $sp_\lambda$ and regard the corresponding symplectic Schur function as a polynomial in $(\bm p)$.
\end{rmk}

The importance of Schur functions is that they encode information on irreducible representations of symmetric groups,
as well as the polynomial representations of general linear groups (see \cite{Mac,S99}).
So do symplectic Schur functions.
They encode information on the irreducible polynomial representations of even symplectic groups (see \cite{Sun90b}).
For example,
the following special evaluation of the symplectic Schur functions
(see \cite{CS12,SK79})
\begin{align}\label{eqn:eval spn}
	sp_{\lambda,2n}(1^n,1^n)
	=\prod_{u\in\lambda}
	\frac{2n+c_{sp}(u)}{h(u)}
\end{align}
is related to the dimension of corresponding representations of symplectic groups,
where $l(\lambda)\leq n$.
It is also related to Amdeberhan's conjectures.
Thus,
under the $p_k(\tilde{\bm x})$ coordinates,
equation \eqref{eqn:eval spn} could be rewritten as
\begin{align}\label{eqn:eval sp}
	sp_{\lambda}(\tilde{\bm x})|_{p_k(\tilde{\bm x})=2n, \forall k}
	=sp_{\lambda}(\bm p)|_{p_k=2n, \forall k}
	=\prod_{u\in\lambda}
	\frac{2n+c_{sp}(u)}{h(u)}.
\end{align}
Since both sides of the above equation could be treated as polynomials in $n$ of finite orders,
and since the above equation holds for any positive integer $n\geq l(\lambda)$,
it should still be valid when replacing $n$ by any complex number $t$.

\subsection{Inner product and the vacuum expectation value}
There is a standard inner product on the ring $\mathbb{C}[p_1,p_2,...]$ by letting $(s_\lambda(\bm p))_{\lambda\in\mathcal{P}}$ be a standard normalized basis (see \cite{Mac} for more details).

Below,
we introduce the notation of vacuum expectation value,
which is more convenient for our purposes (see, for example, \cite{DJM}).
To avoid the misuse of symbols,
we employ the $(\bm q)=(q_k)$ coordinates and vectors that appear in vacuum expectation value are viewed as elements in $\mathbb{C}[\bm q]$.
First,
denote by $|\lambda\rangle$ the Schur function $s_\lambda(\bm q)$ when regarding it as a vector in the vector space $\mathbb{C}[\bm q]$,
and $\langle\lambda|$ is its dual.
The vacuum $|0\rangle$ denotes the vector corresponding to the Schur function labeled by the empty partition,
which is the constant function 1,
and $\langle0|$ is its dual.
Then the vacuum expectation value is of the following form
\begin{align*}
	\langle\lambda|
	A
	|\mu\rangle
	:=(|\lambda\rangle,A\cdot|\mu\rangle)
	=(A^*\cdot|\lambda\rangle,|\mu\rangle),
\end{align*}
where $A$ is an operator acting on the ring $\mathbb{C}[\bm q]$,
and $(\cdot,\cdot)$ is the standard inner product.

\subsection{Vertex operator}
In this subsection,
we review the concept of vertex operators.
It is an important tool in studying the Schur functions and related symmetric functions.
We mainly follow the notations in \cite{DJM}.

Vertex operators are generating functions of differential operators on the ring $\mathbb{C}[\bm q]$.
First, introduce the Heisenberg operators
\begin{align*}
	\alpha_n:=\begin{cases}
		n\frac{\partial}{\partial q_n}, &n>0,\\
		q_{-n} \cdot, &n<0,
	\end{cases}
\end{align*}
where $q_{-n} \cdot$ means the operation of multiplying by $q_{-n}$.
Then, the vertex operators can be defined directly by
\begin{align*}
	\Gamma_\pm(z)=\exp\Big(\sum_{n=1}^\infty \frac{z^n \alpha_{\pm n}}{n}\Big).
\end{align*}
The dual operator of $\Gamma_\pm(z)$ is $\Gamma_\mp(z)$.
By the standard Heisenberg relations $[\alpha_m,\alpha_n]=m\delta_{m+n,0}$,
one can check that
\begin{align}\label{eqn:gamma+- comm relation}
	\Gamma_+(z)\Gamma_-(w)
	=\frac{1}{1-zw}
	\cdot \Gamma_-(w)\Gamma_+(z).
\end{align}

The so called grading operator on the ring $\mathbb{C}[\bm q]$ is
\begin{align*}
	L_0:=\sum_{k=1}^\infty kq_k\frac{\partial}{\partial q_k}.
\end{align*}
It is the homogeneous operator if we regard $\mathbb{C}[\bm q]$ as a graded ring by assigning the degree of $q_k$ to be $k$.
The commutation relations between $L_0$ and $\Gamma_\pm(z)$ are
\begin{align}\label{eqn:L0 gamma+_}
	q^{L_0} \ \Gamma_-(z)
	=\Gamma_-(qz) q^{L_0}
	\quad\text{and}\quad
	q^{L_0} \ \Gamma_+(z)
	=\Gamma_+(q^{-1}z) q^{L_0}.
\end{align}

For a sequence of variables $(\bm p)=(p_1,p_2,...)$,
the corresponding vertex operators are defined by
\begin{align*}
	\Gamma_\pm(\bm p)=\exp\Big(\sum_{n=1}^\infty \frac{p_n \alpha_{\pm n}}{n}\Big).
\end{align*}
When regarding $(\bm p)$ as the power sum coordinates $\big(\bm p(\bm y)\big)=\big(p_k(\bm y)\big)$, in the variables $(\bm y)=(y_1,y_2,...)$,
one has
\begin{align}\label{eqn:def gammapy}
	\Gamma_\pm\big(\bm p(\bm y)\big)
	=\prod_{j=1}^\infty \Gamma_\pm(y_j).
\end{align}
The relation between vertex operators and Schur functions is fundamental (see, for example, equation (A.15) in \cite{O01}) 
\begin{align}\label{eqn:skew schu as vev}
	s_{\lambda/\mu}\big(\bm p(\bm y)\big)
	=\langle\lambda|\Gamma_-\big(\bm p(\bm y)\big)|\mu\rangle
	=\langle\mu|\Gamma_+\big(\bm p(\bm y)\big)|\lambda\rangle.
\end{align}
The special cases of the above equations give
\begin{align}\label{eqn:gamma preserve vacuum}
	\Gamma_+(z)|0\rangle=|0\rangle
	\quad\text{and}\quad
	\langle0|\Gamma_-(z)=\langle0|,
\end{align}
which could easily be proven from their definitions.

\section{A partition function consisting of symplectic Schur functions}
\label{sec: proof single}
In this section,
by studying a partition function consisting of symplectic Schur functions,
we prove Theorem \ref{thm:main csp} and Corollary \ref{cor:main co}.

\subsection{A partition function of symplectic Schur functions}
\begin{thm}\label{thm:sum sp as p_k}
	We have
	\begin{align}\label{eqn:sum sp as p_k}
		\begin{split}
			\sum_{\lambda\in\mathcal{P}}
			q^{|\lambda|}
			\cdot sp_{\lambda}&(\bm p)
			=\prod_{n=1}^\infty \Bigg(\big(1+(-q^2)^n\big)
			\cdot \exp\Big(\sum_{k=1}^\infty
			\frac{(q^{(4n-3)k}-q^{(4n-1)k})p_k-q^{8nk}p_{2k}}
			{k}\Big)\\
			&\cdot \exp\Big(\sum_{k=1}^\infty
			\frac{(q^{(8n-6)k}-q^{(8n-4)k}-q^{8nk})
				\big(p_k^2-p_{2k}\big)
				+q^{(8n-2)k}
				\big(p_k^2+p_{2k}\big)} {2k}\Big)\Bigg).
		\end{split}
	\end{align}
	Both sides of the above equation could be regarded as elements in the ring $\mathbb{C}[\![q]\!][\bm p]$.
\end{thm}
Theorem \ref{thm:sum sp as p_k} has the following equivalent formalism expressed in terms of the variables $(\bm y)=(y_i)$.
\begin{lem}\label{cor:sum sp as p_ky}
	Write $p_k$ as $p_k(\bm y)$, the $k$-power sum of variables $(\bm y)$.
	Then, Theorem \ref{thm:sum sp as p_k} is equivalent to the following
	\begin{align}\label{eqn:sum sp as y_j}
		\begin{split}
			\sum_{\lambda\in\mathcal{P}}
			q^{|\lambda|}
			\cdot &sp_{\lambda}(\bm p(\bm y))
			=\prod_{n=1}^\infty \big(1+(-q^2)^n\big)\\
			&\ \ \cdot \prod_{n=1}^\infty
			\frac{\prod_{1\leq i<j}
				(1-q^{8n-4}y_iy_j)(1-q^{8n}y_iy_j)
				\cdot\prod_{i=1}^\infty
				(1-q^{4n-1}y_i)(1-q^{8n}y_i^2)}
			{\prod_{i=1}^\infty(1-q^{4n-3}y_i)
				\cdot\prod_{1\leq i<j}(1-q^{8n-6}y_iy_j)
				\cdot\prod_{1\leq i\leq j}(1-q^{8n-2}y_iy_j)}.
		\end{split}
	\end{align}
\end{lem}
{\bf Proof:}
From the theory of symmetric functions and the definition of $p_k(\bm y)$,
we have
\begin{align}\label{eqn:symm y as p}
	\begin{split}
		\prod_{i=1}^\infty (1-ay_i)
		=&\exp\Big(-\sum_{k=1}^\infty \frac{a^k}{k}p_k(\bm y)\Big),\\
		\prod_{i=1}^\infty (1-ay_i^2)
		=&\exp\Big(-\sum_{k=1}^\infty \frac{a^k}{k}p_{2k}(\bm y)\Big),\\
		\prod_{1\leq i<j} (1-ay_iy_j)
		=&\exp\Big(-\sum_{k=1}^\infty
		\frac{a^k}{2k}\big(p_k(\bm y)\big)^2-p_{2k}(\bm y)\big)\Big),\\
		\prod_{1\leq i\leq j} (1-ay_iy_j)
		=&\exp\Big(-\sum_{k=1}^\infty
		\frac{a^k}{2k}\big(p_k(\bm y)\big)^2+p_{2k}(\bm y)\big)\Big).
	\end{split}
\end{align}
Then the right hand side of equation \eqref{eqn:sum sp as y_j} is a direct rewriting of equation \eqref{eqn:sum sp as p_k} in terms of the variables $(\bm y)=(y_i)$.
$\Box$

\subsection{Proof of Theorem \ref{thm:sum sp as p_k}}
In this subsection,
we prove Theorem \ref{thm:sum sp as p_k} by showing equation \eqref{eqn:sum sp as y_j}.

For two given operators $A, B$,
if we already know their commutation relation,
then the following special case of the well-known Baker--Campbell--Hausdorff formula (see, for example, Lemma 5.3 in \cite{Mil})
\begin{align}\label{eqn:bch}
	e^A B e^{-A}
	=\sum_{n=0}^\infty \frac{1}{n!} ad_A^n(B),
\end{align}
where $ad_A^0(B):=B, ad_A^n(B)=[A, ad_A^{n-1}(B)]$,
tells us how to compute the commutation relation involving their exponents.
Fortunately,
even if there is an infinite sum in this formula,
in all of our cases,
we have that $ad_A^2(B)$ is a scalar operator then $ad_A^n(B)=0$ for all $n>2$.
Moreover, if $[A,B]$ is already a scalar operator,
then the Zassenhaus formula reads
\begin{align}\label{eqn:z}
	e^{A+B}=e^A e^B e^{-[A,B]/2}.
\end{align}
These two formulas are the main tools in the following Lemma \ref{lem: gamma+- G}.

To prove equation \eqref{eqn:sum sp as y_j},
we first need to introduce the following operator
\begin{align}\label{eqn:def G}
	G:=-\sum_{n=1}^\infty \frac{1}{2n}(\alpha_{-n}^2-\alpha_{-2n}).
\end{align}
The motivation to consider such an operator comes from equations \eqref{eqn:sla=y} and \eqref{eqn:y=p},
which will be exhibited later.
Below,
we first study the commutation relations with respective to $\exp(G)$.
\begin{lem}\label{lem: gamma+- G}
	The operators $\Gamma_-(z)$ and $\exp(G)$
	commute with each other.
	For $\Gamma_+(z)$,
	we have the following commutation relations
	\begin{align}
		\Gamma_+(z) \exp(G)
		=&\exp(G)
		\Gamma_-(z)^{-1}\Gamma_+(z), \label{eqn: gamma+ sla}\\
		\Gamma_+(z)^{-1} \exp(G)
		=&(1-z^2)
		\cdot \exp(G)
		\Gamma_-(z)\Gamma_+(z)^{-1}. \label{eqn: gamma+-1 sla}
	\end{align}
\end{lem}
{\bf Proof:}
First,
since both of the operators $\Gamma_-(z)$ and $\exp\Big(\sum_{n=1}
\frac{1}{n}(\frac{1}{2}\alpha_{-n}^2-\frac{1}{2}\alpha_{-2n})\Big)$ consist only of $\alpha_{-n}$ for $n>0$,
they commute with each other.
Next,
one can use the Heisenberg relations $[\alpha_m,\alpha_n]=m\delta_{m+n,0}$ to obtain
\begin{align*}
	\big[\sum_{n=1}^\infty \frac{z^n \alpha_{n}}{n},
	\sum_{n=1}^{\infty}
	\frac{1}{2n}\alpha_{-2n}\big]
	=\sum_{n=1}^{\infty}
	\frac{z^{2n}}{2n}
	=-\frac{1}{2}\ln(1-z^2).
\end{align*}
Then,
one can directly apply the Baker--Campbell--Hausdorff formula \eqref{eqn:bch} to obtain the following commutation relation
\begin{align}\label{eqn:gamma+1}
	\Gamma_+(z)\exp\Big(\sum_{n=1}^{\infty}
	\frac{1}{2n}\alpha_{-2n}\Big)
	=(1-z^2)^{-1/2}
	\cdot\exp\Big(\sum_{n=1}^{\infty}
	\frac{1}{2n}\alpha_{-2n}\Big)
	\Gamma_+(z).
\end{align}
By the similar method,
one can derive
\begin{align}\label{eqn:gamma+2}
	\Gamma_+(z)\exp\Big(-\sum_{n=1}^{\infty}
	\frac{1}{2n}\alpha_{-n}^2\Big)
	=(1-z^2)^{1/2}
	\cdot\exp\Big(-\sum_{n=1}^{\infty}
	\frac{1}{2n}\alpha_{-n}^2\Big)
	\Gamma_-(z)^{-1}\Gamma_+(z)
\end{align}
using the Baker--Campbell--Hausdorff formula \eqref{eqn:bch}.
It additionally demands the following 
\begin{align}
	\exp\Big(\sum_{n=1}^{\infty}
	\frac{z^n}{n}(\alpha_n-\alpha_{-n})\Big)
	=(1-z^2)^{-1/2}
	\cdot\Gamma_-(z)^{-1}\Gamma_+(z),
\end{align}
which is provable by the Zassenhaus formula \eqref{eqn:z}.
Then,
equation \eqref{eqn: gamma+ sla} is a consequence of the above two commutation relations \eqref{eqn:gamma+1} and \eqref{eqn:gamma+2}.
Equation \eqref{eqn: gamma+-1 sla} can be shown by multiplying both sides of equation \eqref{eqn: gamma+ sla} by $\Gamma_+(z)^{-1}$.
$\Box$

\begin{lem}\label{lem: gamma+ G}
	For vertex operator $\Gamma_+\big(\bm p(\bm y)\big)$, defined in equation \eqref{eqn:def gammapy},
	we have the following commutation relations
	\begin{align}
		\Gamma_+\big(\bm p(\bm y)\big) \exp(G)
		=&\prod_{1\leq i<j}(1-y_iy_j)
		\cdot\Gamma_-\big(\bm p(\bm y)\big)^{-1}
		\exp(G)
		\Gamma_+\big(\bm p(\bm y)\big), \label{eqn: gamma+ G}\\
		\Gamma_+\big(\bm p(\bm y)\big)^{-1} \exp(G)
		=&\prod_{i=1}^\infty(1-y_i^2)
		\cdot\prod_{1\leq i<j}(1-y_iy_j)
		\cdot\Gamma_-\big(\bm p(\bm y)\big)
		\exp(G)
		\Gamma_+\big(\bm p(\bm y)\big)^{-1}. \label{eqn: gamma+-1 G}
	\end{align}
\end{lem}
{\bf Proof:}
Just recalling that
$\Gamma_+\big(\bm p(\bm y)\big)=\prod_{j=i}^\infty \Gamma(y_i)$
and using Lemma \ref{lem: gamma+- G},
we have
\begin{align*}
	\Gamma_+\big(\bm p(\bm y)\big) \exp(G)
	=&\exp(G)
	\prod_{i=1}^\infty\big(\Gamma_-(y_i)^{-1}\Gamma_+(y_i)\big),\\
	\Gamma_+\big(\bm p(\bm y)\big)^{-1} \exp(G)
	=&\prod_{i=1}^\infty(1-y_i^2)
	\cdot\exp(G)
	\prod_{i=1}^\infty\big(\Gamma_-(y_i)\Gamma_+(y_i)^{-1}\big).
\end{align*}
Then, the assertions follow from the commutation relation \eqref{eqn:gamma+- comm relation}.
$\Box$

Motivated from the I.5.Example 4 in \cite{Mac} (see also equation (5.2) in \cite{Y23} for the following equivalent form),
\begin{align*}
	\sum_{\lambda\in\mathcal{P}} s_\lambda(\bm p)=
	\exp\Big(\sum_{n=1}^\infty \frac{1}{n}\big(p_n+\frac{1}{2}p_n^2-\frac{1}{2}p_{2n}\big)\Big),
\end{align*}
it is beneficial to introduce another useful operator as
\begin{align}\label{eqn:def F}
	F:=\sum_{n=1}^\infty
	\frac{1}{n}(\alpha_{-n}+\frac{1}{2}\alpha_{-n}^2-\frac{1}{2}\alpha_{-2n}).
\end{align}
Denote its dual by
$F^*:=\sum_{n=1}^\infty
\frac{1}{n}(\alpha_{n}+\frac{1}{2}\alpha_{n}^2-\frac{1}{2}\alpha_{2n})$.
Then we have

\begin{lem}\label{lem: gamma+- F}
	The operators $\Gamma_+(\bm p(\bm y))$ and $\exp(F^*)$ commute with each other.
	For $\Gamma_-(\bm p(\bm y))$,
	we have
	\begin{align}
		\exp(F^*) \Gamma_-\big(\bm p(\bm y)\big)
		=&\frac{1}{\prod\limits_{i=1}^\infty(1-y_i)
			\cdot\prod\limits_{1\leq i<j}(1-y_iy_j)}
		\cdot \Gamma_-\big(\bm p(\bm y)\big)
		\exp(F^*)
		\Gamma_+\big(\bm p(\bm y)\big), \label{eqn: sla gamma-}\\
		\exp(F^*) \Gamma_-\big(\bm p(\bm y)\big)^{-1}
		=&\frac{\prod_{i=1}^\infty (1-y_i)}
		{\prod_{1\leq i\leq j}(1-y_iy_j)}
		\cdot \Gamma_-\big(\bm p(\bm y)\big)^{-1} \exp(F^*)\Gamma_+\big(\bm p(\bm y)\big)^{-1}. \label{eqn: sla gamma-^-1}
	\end{align}
\end{lem}
{\bf Proof:} The author proved (see  Lemma 5.2 in \cite{Y23}) that
\begin{align*}
	\Gamma_+(z)
	\exp(F)
	=\frac{1}{1-z}
	\cdot\exp(F)
	\Gamma_-(z)\Gamma_+(z).
\end{align*}
By using the dual of the above equation, one can show
\begin{align*}
	\exp(F^*) \Gamma_-\big(\bm p(\bm y)\big)
	=&\frac{1}{\prod_{i=1}^\infty(1-y_i)}
	\cdot 
	\prod_{i=1}^\infty \big(\Gamma_-(y_i)\Gamma_+(y_i)\big)
	\cdot \exp(F^*).
\end{align*}
Thus, equation \eqref{eqn: sla gamma-} follows from the above equation and the commutation relation \eqref{eqn:gamma+- comm relation}.
To prove the relation in \eqref{eqn: sla gamma-^-1},
one only needs to multiply both sides of equation \eqref{eqn: sla gamma-} by $\Gamma\big(\bm p(\bm y)\big)^{-1}$.
$\Box$
\ \\

\indent With the above commutation relations at hand, we are prepared to prove Theorem \ref{thm:sum sp as p_k}.

{\bf Proof of Theorem \ref{thm:sum sp as p_k}:}
Denote by $\mathcal{P}'$ the set of partitions whose Frobenius notations are of the form $(m_1,...,m_r|m_1+1,...,m_r+1)$ for some decreasing non-negative integers $m_1,...,m_r$.
Actually,
after taking the conjugation,
this kind of partition gives the doubled distinct partition studied in \cite{GKS,Lit}.
For convenience, the empty partition $\emptyset$ is considered in $\mathcal{P}'$.
The following result was proved by Sundaram in her thesis
(see Theorem 12.10 in \cite{Sun86})
\begin{align}\label{eqn:sp as s}
	sp_\lambda(\bm p)
	=\sum_{\mu\subseteq\lambda, \mu\in\mathcal{P}}
	s_{\mu}(\bm p)
	\cdot
	\big(
	\sum_{\beta\in\mathcal{P}'}(-1)^{|\beta|/2}c^{\lambda}_{\mu,\beta}\big),
	\ \ \ \forall\lambda\in\mathcal{P},
\end{align}
which is motivated by a classical result representing Schur functions in terms of symplectic Schur functions obtained by Littlewood (see, for example, \cite{Lit}).
Sundaram provided the equation \eqref{eqn:sp as s} by regarding both sides of the equation as formal functions in variables $(x_1,x_1^{-1},...,x_n,x_n^{-1})$ when $l(\lambda)\leq n$.
Here,
equation \eqref{eqn:sp as s} is an equality in the ring $\mathbb{C}[\bm p]$.

Then from the definition \eqref{eqn:def skew schur} of skew Schur functions,
we have
\begin{align}\label{eqn:sum sp as skew s}
	\begin{split}
		\sum_{\lambda\in\mathcal{P}}
		q^{|\lambda|}
		\cdot sp_{\lambda}(\bm p)
		=&\sum_{\lambda\in\mathcal{P} \atop \beta\in\mathcal{P}'}
		q^{|\lambda|} (-1)^{|\beta|/2} s_{\lambda/\beta}(\bm p)\\
		=&\sum_{\lambda\in\mathcal{P} \atop \beta\in\mathcal{P}'}
		\langle\lambda|q^{L_0}
		\cdot \Gamma_-(\bm p)
		\cdot (-1)^{L_0/2} |\beta\rangle,
	\end{split}
\end{align}
where the second equality sign is an application of formula \eqref{eqn:skew schu as vev},
which expresses the skew Schur function in terms of the vacuum expectation value involving vertex operators.

The following equation is useful 
(see I.5.Example 9 (a) in \cite{Mac}),
\begin{align}\label{eqn:sla=y}
	\sum_{\lambda\in\mathcal{P}'} (-1)^{|\lambda|/2}
	s_{\lambda}\big(\bm p(\bm y)\big)
	=\prod_{i<j} (1-y_i y_j).
\end{align}
When expressed in terms of the power sum coordinates $\big(\bm p(\bm y)\big)$,
one has
\begin{align}\label{eqn:y=p}
	\prod_{i<j} (1-y_i y_j)
	=\exp\Big(-\sum_{n=1}^\infty \frac{1}{2n}\big(p_n(\bm y)^2-p_{2n}(\bm y)\big)\Big).
\end{align}
The operator $G$, defined in equation \eqref{eqn:def G}, is indeed the operator by multiplying the function above.
Moreover,
together with the homogeneous operator $L_0$,
the above two equations \eqref{eqn:sla=y} and \eqref{eqn:y=p} can be rewritten as
\begin{align}\label{eqn:sum slambda' G}
	\sum_{\lambda\in\mathcal{P}'} (-1)^{L_0/2}
	|\lambda\rangle=
	\exp(G)
	\cdot |0\rangle.
\end{align}
For the operator $F$,
we have the following similar formulas
(see I.5.Example 4 in \cite{Mac},
and see also equation (5.4) in \cite{Y23}),
\begin{align}\label{eqn:sum sla |la>}
	\sum_{\lambda\in\mathcal{P}} |\lambda\rangle=
	\exp(F)
	\cdot |0\rangle,
	\quad\quad
	\sum_{\lambda\in\mathcal{P}} \langle\lambda|=
	\exp(F^*)
	\cdot \langle0|.
\end{align}

Thus, equation \eqref{eqn:sum sp as skew s} can be expressed as
\begin{align}\label{eqn:sum sp as vev}
	\sum_{\lambda\in\mathcal{P}}
	q^{|\lambda|}
	\cdot sp_{\lambda}(\bm p)
	=\langle0|\exp(F^*)
	\cdot q^{L_0} \cdot \Gamma_-(\bm p)
	\cdot \exp(G)|0\rangle.
\end{align}
From now on, and for convenience,
we replace the variables $(\bm p)$ by symmetric functions $\big(\bm p(\bm y)\big)$.
One can first use equation \eqref{eqn: sla gamma-} to change the orders of the operators $\Gamma_-\big(\bm p(\bm y)\big)$ and $\exp(F^*)$.
The result is that,
$\Gamma_-\big(\bm p(\bm y)\big)$ becomes $\Gamma_+\big(\bm p(q^2\bm y)\big)$ since $\langle0|\Gamma_-(\cdot)=\langle0|$.
Then one can use equation \eqref{eqn: gamma+ G} to deal with the commutation relation between operators $\Gamma_+\big(\bm p(q^2\bm y)\big)$ and $\exp(G)$.
More precisely, the left hand side of equation \eqref{eqn:sum sp as vev} is equal to
\begin{align*}
	&\langle0|\exp(F^*)\cdot q^{L_0}\cdot \Gamma_+\big(\bm p(q^2\bm y)\big)
	\cdot\exp(G)|0\rangle
	\cdot\frac{1}
	{\prod_{i=1}^\infty(1-qy_i)
		\cdot\prod_{1\leq i<j}(1-q^2y_iy_j)}\\
	=&\langle0|\exp(F^*)\cdot q^{L_0}\cdot \Gamma_-\big(\bm p(q^2\bm y)\big)^{-1}
	\cdot\exp(G)|0\rangle
	\cdot \frac{\prod_{1\leq i<j}(1-q^4y_iy_j)}{\prod_{i=1}^\infty(1-qy_i)
		\cdot\prod_{1\leq i<j}(1-q^2y_iy_j)},
\end{align*}
where we have used the fact that $\Gamma_+(\cdot)|0\rangle=|0\rangle$.
Similarly,
one can use equations \eqref{eqn: sla gamma-^-1} and \eqref{eqn: gamma+-1 G} to change the orders of the operators $\Gamma_\pm(\cdot)^{-1}$, $\exp(F^*)$ and $\exp(G)$.
The result is that,
the left hand side of equation \eqref{eqn:sum sp as vev} is equal to
\begin{align*}
	&\langle0|\exp(F^*) q^{L_0} \Gamma_+\big(\bm p(q^4\bm y)\big)^{-1}
	\exp(G)|0\rangle
	\cdot \frac{\prod_{1\leq i<j}(1-q^4y_iy_j)
		\cdot\prod_{i=1}^\infty(1-q^3y_i)}
	{\prod\limits_{i=1}^\infty(1-qy_i)
		\cdot\prod\limits_{1\leq i<j}(1-q^2y_iy_j)
		\cdot\prod\limits_{1\leq i\leq j}(1-q^6y_iy_j)}\\
	=&\langle0|\exp(F^*) q^{L_0} \Gamma_-\big(\bm p(q^4\bm y)\big)
	\exp(G)|0\rangle\\
	&\quad\quad\cdot
	\frac{\prod_{1\leq i<j}(1-q^4y_iy_j)
		\cdot\prod_{i=1}^\infty(1-q^3y_i)
		\cdot \prod_{i=1}^\infty (1-q^8y_i^2)
		\cdot \prod_{1\leq i<j}(1-q^8y_iy_j)}
	{\prod_{i=1}^\infty(1-qy_i)
		\cdot\prod_{1\leq i<j}(1-q^2y_iy_j)
		\cdot\prod_{1\leq i\leq j}(1-q^6y_iy_j)}.
\end{align*}
This time, we compare the above equation and equation \eqref{eqn:sum sp as vev}.
One can notice that,
for an arbitrary positive integer $M$,
by repeating the preceding process $M$ times,
the following equation can be obtained
\begin{align}\label{eqn:first main last eqn}
	\begin{split}
		\sum_{\lambda\in\mathcal{P}}
		q^{|\lambda|}
		\cdot sp_{\lambda}&\big(\bm p(\bm y)\big)
		=\langle0|\exp(F^*) q^{L_0}
		\Gamma_-\big(\bm p(q^{4M}\bm y)\big)
		\exp(G)|0\rangle\\
		&\cdot \prod_{n=1}^M
		\frac{\prod\limits_{1\leq i<j}
			(1-q^{8n-4}y_iy_j)(1-q^{8n}y_iy_j)
			\cdot\prod\limits_{i=1}^\infty
			(1-q^{4n-1}y_i)(1-q^{8n}y_i^2)}
		{\prod\limits_{i=1}^\infty(1-q^{4n-3}y_i)
			\cdot\prod\limits_{1\leq i<j}(1-q^{8n-6}y_iy_j)
			\cdot\prod\limits_{1\leq i\leq j}(1-q^{8n-2}y_iy_j)}.
	\end{split}
\end{align}
Since both sides of the above equation are regarded as elements in the ring $\Lambda(\bm y)[\![q]\!]$,
we have,
as operators on the ring $\Lambda(\bm y)[\![q]\!]$,
\[\Gamma_-\big(\bm p(q^{4M}\bm y)\big)\rightarrow \text{Id},
\ \ \ \text{when}\ \ \ M\rightarrow+\infty,\]
where $\text{Id}$ denotes the identity operator.
Thus,
the theorem is proved by taking the limit $M\rightarrow+\infty$ in equation \eqref{eqn:first main last eqn} and by applying the following combinatorial result
\begin{align*}
	\langle0|\exp(F^*)\cdot q^{L_0}
	\cdot\exp(G)|0\rangle
	=\sum_{\lambda\in\mathcal{P}'}
	(-1)^{|\lambda|/2} q^{|\lambda|}
	=\prod_{n=1}^\infty \big(1+(-q^2)^n\big),
\end{align*}
where the first equality follows from equations \eqref{eqn:sum slambda' G} and \eqref{eqn:sum sla |la>},
while the second equality follows from the fact that the size of the partition $(m_1,...,m_r|m_1+1,...,m_r+1)\in\mathcal{P}'$
is equal to $2\sum_{i=1}^r(m_i+1)$.
$\Box$

\subsection{Hook-symplectic/orthogonal contents formulas}
In this subsection,
we venture on proving the results in Theorem \ref{thm:main csp} and Corollary \ref{cor:main co}.
\ \\

\indent{\bf Proof of Theorem \ref{thm:main csp}:}
What we need is to consider the specialization $p_k=t, \forall k\geq1$ of equation \eqref{eqn:sum sp as p_k}.
Equivalently,
for an arbitrary positive integer $n_1$,
setting
\[y_i=1, 1\leq i \leq n_1,
\ \ \ \text{and}\ \ \ y_i=0, i>n_1\]
in equation \eqref{eqn:sum sp as y_j},
by the evaluation formula \eqref{eqn:eval sp} for the symplectic Schur functions,
we obtain
\begin{align}\label{eqn:first main as n1}
	\begin{split}
		\sum_{\lambda\in\mathcal{P}}q^{|\lambda|}
		\prod_{u\in\lambda}
		&\frac{n_1+c_{sp}(u)}{h(u)}\\
		&=\prod_{n=1}^\infty
		\Bigg(\frac{1-q^{4n-2}}{1-q^{8n-4}}
		\cdot\frac{(1-q^{8n-4})^{\binom{n_1}{2}}
			(1-q^{4n-1})^{n_1}
			(1-q^{8n})^{\binom{n_1}{2}+n_1}}
		{(1-q^{4n-3})^{n_1}
			(1-q^{8n-6})^{\binom{n_1}{2}}
			(1-q^{8n-2})^{\binom{n_1+1}{2}}}\Bigg).
	\end{split}
\end{align}
Thus,
we indeed have proved Theorem \ref{thm:main csp} when $t$ is a positive integer.
This is enough.
Actually,
one can expand both sides of the above equation \eqref{eqn:first main as n1} in the formal variable $q$.
Then the coefficients on the two sides could be regarded as polynomials in $n_1$ of finite orders.
Since we proved the equivalence of these two polynomials evaluated at infinite many positive integers $n_1$,
they should be equal to each other as polynomials.
That is to say,
the above equation holds when replacing $n_1$ by any complex number $t$.
After doing that,
the result is equivalent to equation \eqref{eqn:main csp}.
$\Box$
\ \\

\indent{\bf Proof of Corollary \ref{cor:main co}:}
Starting with the definitions \eqref{eqn:def csp} and \eqref{eqn:def co} of symplectic and orthogonal contents,
an easy combinatorial fact arises as
\begin{align*}
	h^\lambda(i,j)=h^{\lambda^t}(j,i)
	\text{\ \ \ \ and\ \ \ \ }
	c^\lambda_{sp}(i,j)
	=-c^{\lambda^t}_{O}(j,i).
\end{align*}
As a consequence,
the weight functions on the left hand side of equations \eqref{eqn:main csp} and \eqref{eqn:main co} are related to each other by
\begin{align*}
	\prod_{u\in\lambda}
	\frac{t+c_{O}(u)}{h(u)}
	=(-1)^{|\lambda|}\prod_{u\in\lambda^t}
	\frac{-t+c_{sp}(u)}{h(u)}
\end{align*}
for an arbitrary partition $\lambda$ and its conjugation $\lambda^t$.
Since taking conjugation is an involution in the set of partitions,
we have
\begin{align*}
	\sum_{\lambda\in\mathcal{P}}q^{|\lambda|}
	\prod_{u\in\lambda}
	\frac{t+c_{O}(u)}{h(u)}
	=\big(\text{r.h.s.\ of\ equation\ \eqref{eqn:main csp}}\big)
	|_{t\mapsto-t, q\mapsto-q},
\end{align*}
which proves the corollary just by noticing
\begin{align*}
	\prod_{n=1}^\infty
	\Big(\frac{1-q^{4n-1}}{1-q^{4n-3}}\Big)^t
	|_{t\mapsto-t, q\mapsto-q}
	=\prod_{n=1}^\infty
	\Big(\frac{1-q^{4n-1}}{1-q^{4n-3}}
	\cdot \frac{1-q^{8n-6}}{1-q^{8n-2}}\Big)^t.
\end{align*}
$\Box$

\section{A Cauchy-identity type formula for symplectic Schur functions}
\label{sec: proof double}

\subsection{Partition function of double symplectic Schur functions}
In this subsection,
we obtain a formula for the partition function consisting of double symplectic Schur functions as in Theorem \ref{thm:sum spsp as p_k}.

Denote the dual operator of $G$ by
\begin{align*}
	G^*:=-\sum_{n=1}^\infty \frac{1}{2n}(\alpha_{n}^2-\alpha_{2n}).
\end{align*}
We first need the following commutation relation.
\begin{cor}\label{cor: gamma- G^*}
	The operators $\Gamma_+(z)$ and $\exp(G^*)$ commute with each other.
	For the vertex operator $\Gamma_-\big(\bm p(\bm y)\big)$,
	we have
	\begin{align}
		\exp(G^*) \Gamma_-\big(\bm p(\bm y)\big)
		=&\prod_{1\leq i<j}(1-y_iy_j)
		\cdot\Gamma_-\big(\bm p(\bm y)\big)
		\exp(G^*)
		\Gamma_+\big(\bm p(\bm y)\big)^{-1} \label{eqn: gamma- G*},\\
		\exp(G^*)
		\Gamma_-\big(\bm p(\bm y)\big)^{-1}
		=&\prod_{i=1}^\infty(1-y_i^2)
		\cdot\prod_{1\leq i<j}(1-y_iy_j)
		\cdot\Gamma_-\big(\bm p(\bm y)\big)^{-1}
		\exp(G^*)
		\Gamma_+\big(\bm p(\bm y)\big). \label{eqn: gamma--1 G*}
	\end{align}
\end{cor}
{\bf Proof:}
Take the dual of Lemma \ref{lem: gamma+ G}.
$\Box$

\begin{thm}\label{thm:sum spsp as p_k}
	We have
	\begin{align*}
		\begin{split}
			\sum_{\lambda\in\mathcal{P}}
			q^{|\lambda|}
			\cdot sp_{\lambda}(\bm p)&sp_{\lambda}(\bm p')
			=\prod_{n=1}^\infty
			\Bigg(
			\exp\Big(\sum_{k=1}^\infty
			\frac{(q^{(4n-3)k}+q^{(4n-1)k})p_kp_k'-q^{4nk}(p_{2k}+p_{2k}')}{k}\Big)\\
			&\cdot\exp\Big(-\sum_{k=1}^\infty
			\frac{(q^{(4n-2)k}+q^{4nk})(p_k^2-p_{2k}+p_k'^2-p'_{2k})}{2k}\Big)
			\cdot \frac{1}{1-q^{4n-2}}\Bigg).
		\end{split}
	\end{align*}
\end{thm}
Use the variables $(\bm y)=(y_i)$ and write $p_k$ as $p_k(\bm y)$, the $k$-th power sum of the variables $(\bm y)$.
Then, we have
\begin{lem}\label{cor:sum spsp as p_k}
	The Theorem \ref{thm:sum spsp as p_k} is equivalent to the following
	\begin{align}\label{eqn:sum spsp as p_ky}
		\begin{split}
			\sum_{\lambda\in\mathcal{P}}
			q^{|\lambda|}
			\cdot sp_{\lambda}\big(\bm p&(\bm y)\big)sp_{\lambda}\big(\bm p'(\bm y')\big)
			=\prod_{n=1}^\infty
			\Bigg( \prod\limits_{i=1}^\infty (1-q^{4n}y_i^2)(1-q^{4n}y_i'^{2})\\
			&\cdot\frac{\prod\limits_{1\leq i<j} (1-q^{4n-2}y_i'y_j')(1-q^{4n-2}y_iy_j)(1-q^{4n}y_iy_j)(1-q^{4n}y_i'y_j')}
			{(1-q^{4n-2})
				\cdot\prod_{i,j=1}^\infty (1-q^{4n-3}y_iy_j') (1-q^{4n-1}y_iy_j')}\Bigg).
		\end{split}
	\end{align}
\end{lem}
{\bf Proof:}
This is similar to the proof of Lemma \ref{cor:sum sp as p_ky}.
One just needs to notice equation \eqref{eqn:symm y as p} and the following fact
\begin{align*}
	\prod_{i,j=1}^\infty (1-ay_iy_j')
	=\exp\Big(-\sum_{k=1}^\infty\frac{a^k}{k}p_k(\bm y)p'_k(\bm y')\Big).
\end{align*}
$\Box$
\ \\

\indent{\bf Proof of Theorem \ref{thm:sum spsp as p_k}:}
Observing that equation \eqref{eqn:sum spsp as p_ky} is equivalent to Theorem \ref{thm:sum spsp as p_k}, we will prove \eqref{eqn:sum spsp as p_ky}.
First,
by equation \eqref{eqn:sp as s},
we have
\begin{align*}
	\sum_{\lambda\in\mathcal{P}}
	q^{|\lambda|}
	\cdot sp_{\lambda}(\bm p)sp_{\lambda}(\bm p')
	=&\sum_{\lambda\in\mathcal{P} \atop \beta,\gamma\in\mathcal{P}'}
	q^{|\lambda|} (-1)^{(|\beta|+|\gamma|)/2} s_{\lambda/\beta}(\bm p)
	s_{\lambda/\gamma}(\bm p')\\
	=&\sum_{\lambda\in\mathcal{P} \atop \beta,\gamma\in\mathcal{P}'}
	\langle\lambda|q^{L_0}
	\cdot \Gamma_-(\bm p)
	\cdot (-1)^{L_0/2} |\beta\rangle
	\cdot \langle\lambda|
	\cdot \Gamma_-(\bm p')
	\cdot (-1)^{L_0/2} |\gamma\rangle.
\end{align*}
By taking the dual of the last term from the right hand side of the above equation,
and using equation \eqref{eqn:sum slambda' G},
the left hand side of the above equation is equal to
\begin{align*}
	&\sum_{\beta,\gamma\in\mathcal{P}'}
	\langle\gamma|(-1)^{L_0/2}
	\cdot \Gamma_+(\bm p')
	\cdot q^{L_0}
	\cdot \Gamma_-(\bm p)
	\cdot (-1)^{L_0/2} |\beta\rangle\\
	=&\langle0|\exp(G^*)
	\cdot \Gamma_+(\bm p')
	\cdot q^{L_0}
	\cdot \Gamma_-(\bm p)
	\cdot \exp(G)|0\rangle.
\end{align*}
The computation of the above vacuum expectation value is similar to the proof of Theorem \ref{thm:sum sp as p_k}. 
First,
we replace the variables $(\bm p)$ and $(\bm p')$ by symmetric functions $\big(\bm p(\bm y)\big)$ and $\big(\bm p'(\bm y')\big)$, respectively.
The commutation relations \eqref{eqn:L0 gamma+_} and \eqref{eqn:gamma+- comm relation} can be used to change the orders of $q^{L_0}$, $\Gamma_+\big(\bm p'(\bm y')\big)$ and $\Gamma_-\big(\bm p(\bm y)\big)$ to obtain
\begin{align*}
	\sum_{\lambda\in\mathcal{P}}
	q^{|\lambda|}
	\cdot sp_{\lambda}\big(\bm p(\bm y)\big)sp_{\lambda}(\bm p'\big(\bm y')\big)
	=&\frac{\langle0|\exp(G^*)
		\cdot \Gamma_-\big(\bm p(q\bm y)\big)
		\cdot q^{L_0}
		\cdot \Gamma_+\big(\bm p'(q\bm y')\big)
		\cdot \exp(G)|0\rangle}
	{\prod_{i,j=1}^\infty (1-qy_iy_j')}.
\end{align*}
Then we can use equations \eqref{eqn: gamma- G*} and \eqref{eqn: gamma+ G} to change the orders of the operators $\Gamma_\pm(\cdot)$, $\exp(G^*)$ and $\exp(G)$.
The result is that,
the left hand side of equation \eqref{eqn:sum spsp as p_ky} is equal to
\begin{align*}
	&\frac{\prod_{1\leq i<j}(1-q^2y_i'y_j')
		\cdot \prod_{1\leq i<j}(1-q^2y_iy_j)}
	{\prod_{i,j=1}^\infty (1-qy_iy_j')}\\
	&\quad\quad\quad\quad\quad\quad\cdot\langle0|\exp(G^*)
	\cdot \Gamma_+\big(\bm p(q\bm y)\big)^{-1}
	\cdot q^{L_0}
	\cdot \Gamma_-\big(\bm p'(q\bm y')\big)^{-1}
	\cdot\exp(G)|0\rangle.
\end{align*}
By analogy,
we can use equations \eqref{eqn: gamma+-1 G} and \eqref{eqn: gamma--1 G*} to deal with the commutation relations between $\Gamma_\pm(\cdot)^{-1}$ and $\exp(G), \exp(G^*)$.
The result is
\begin{align*}
	&\sum_{\lambda\in\mathcal{P}}
	q^{|\lambda|}
	\cdot sp_{\lambda}\big(\bm p(\bm y)\big)
	sp_{\lambda}\big(\bm p'(\bm y')\big)\\
	=&\frac{\prod\limits_{1\leq i<j} (1-q^2y_i'y_j')(1-q^2y_iy_j)(1-q^4y_iy_j)(1-q^4y_i'y_j')
		\cdot \prod\limits_{i=1}^\infty (1-q^4y_i^2)(1-q^4y_i'^{2})}
	{\prod_{i,j=1}^\infty (1-qy_iy_j') (1-q^3y_iy_j')}\\
	&\cdot\langle0|\exp(G^*)
	\cdot \Gamma_+\big(\bm p'(q^2\bm y')\big)
	\cdot q^{L_0}
	\cdot \Gamma_-\big(\bm p(q^2\bm y)\big)
	\cdot\exp(G)|0\rangle.
\end{align*}
By repeating the above process again and again,
one obtains that,
the left hand side of equation \eqref{eqn:sum spsp as p_ky} is equal to
\begin{align*}
	\langle0|\exp(G^*) q^{L_0}
	\exp(G)&|0\rangle
	\cdot \prod_{n=1}^\infty\Bigg(
	\prod\limits_{i=1}^\infty (1-q^{4n}y_i^2)(1-q^{4n}y_i'^{2})\\
	&\cdot
	\frac{\prod\limits_{1\leq i<j} (1-q^{4n-2}y_i'y_j')(1-q^{4n-2}y_iy_j)(1-q^{4n}y_iy_j)(1-q^{4n}y_i'y_j')}
	{\prod_{i,j=1}^\infty (1-q^{4n-3}y_iy_j') (1-q^{4n-1}y_iy_j')}\Bigg).
\end{align*}
The vacuum expectation value in the above equation can be computed in terms of
\begin{align*}
	\langle0|\exp(G^*) q^{L_0}
	\exp(G)|0\rangle
	=\sum_{\lambda\in\mathcal{P}'}
	(-q)^{|\lambda|}
	=\prod_{n=1}^\infty (1+q^{2n})
	=\prod_{n=1}^\infty \frac{1}{1-q^{4n-2}},
\end{align*}
where we have used equation \eqref{eqn:sum slambda' G}.
This completes the proof of the theorem.
$\Box$

\subsection{Proofs of Theorems \ref{thm:main csp^2 cO^2} \ref{thm:main csp*c} and
	Corollary \ref{cor:main csp^2=csp}}
\ \\
\indent{\bf Proof of Theorem \ref{thm:main csp^2 cO^2}:}
For any two positive integers $n_1$ and $n_2$,
taking the special evaluation of equation \eqref{eqn:sum spsp as p_ky} at
$$y_i=y_j'=1, 1\leq i\leq n_1, 1\leq j\leq n_2
\text{\ \ \ \ and\ \ \ \ }
y_i=y_j'=0, i>n_1, j>n_2,$$
we obtain that
\begin{align*}
	&\sum_{\lambda\in\mathcal{P}}
	q^{|\lambda|}
	\cdot \prod_{u\in\lambda}
	\frac{\big(n_1+c_{sp}(u)\big)\big(n_2+c_{sp}(u)\big)}{h^2(u)}\\
	=&\prod_{n=1}^\infty
	\frac{(1-q^{4n})^{n_1}(1-q^{4n})^{n_2}
		\cdot(1-q^{4n-2})^{\binom{n_2}{2}}
		(1-q^{4n-2})^{\binom{n_1}{2}}
		\cdot(1-q^{4n})^{\binom{n_1}{2}}
		(1-q^{4n})^{\binom{n_2}{2}}}
	{(1-q^{4n-2})
		\cdot (1-q^{4n-3})^{n_1n_2} (1-q^{4n-1})^{n_1n_2}}\\
	=&\prod_{n=1}^\infty
	\frac{(1-q^{4n-2})^{\binom{n_1}{2}+\binom{n_2}{2}-1}
		\cdot (1-q^{4n})^{\binom{n_1+1}{2}+\binom{n_2+1}{2}}}
	{(1-q^{4n-3})^{n_1n_2}\cdot (1-q^{4n-1})^{n_1n_2}}.
\end{align*}
Since both sides of the above equation could be regarded as elements in the ring $\mathbb{C}[\![q]\!][n_1,n_2]$ and we already have proved the above equation for all $n_1, n_2 \in\mathbb{Z}_{>0}$, it follows that the equation should still hold when replacing $n_1, n_2$ by arbitrary complex numbers $t_1, t_2$.
Hence, the first equation \eqref{eqn:main spsp} of Theorem~\ref{thm:main csp^2 cO^2} is thus proved.

For the second equation \eqref{eqn:main coco} in this theorem,
it suffices to make the change of variables $(t_1,t_2)\rightarrow(-t_1,-t_2)$ in equation \eqref{eqn:main spsp} and to notice that
\begin{align*}
	\prod_{u\in\lambda}
	\frac{\big(t_1+c_{O}(u)\big)\big(t_2+c_{O}(u)\big)}{h^2(u)}
	=\prod_{u\in\lambda^t}
	\frac{\big(-t_1+c_{sp}(u)\big)\big(-t_2+c_{sp}(u)\big)}{h^2(u)}
\end{align*}
holds for an arbitrary partition $\lambda$ and its conjugation $\lambda^t$.
$\Box$

\begin{ex}\label{exa:conj6.3a}
	Letting $t_1=-t_2=t$ in equation \eqref{eqn:main spsp}, we obtain
	\begin{align*}
		&\sum_{\lambda\in\mathcal{P}}
		q^{|\lambda|}
		\cdot \prod_{u\in\lambda}
		\frac{-t^2+c_{sp}^2(u)}{h^2(u)}\\
		=&\prod_{n=1}^\infty
		\frac{(1-q^{4n-2})^{t^2-1}
			\cdot(1-q^{4n})^{t^2}}
		{(1-q^{4n-3})^{-t^2}
			\cdot (1-q^{4n-1})^{-t^2}}\\
		=&\prod_{n=1}^\infty
		\frac{1}
		{(1-q^{4n-2})
			\cdot(1-q^{n})^{-t^2}},
	\end{align*}
	which is exactly the first result in Conjecture 6.3 (a) of \cite{Amd12}. Similarly, the second result of Conjecture 6.3 (a) in \cite{Amd12} follows from the special case $t_1=-t_2=t$ of equation \eqref{eqn:main coco}.
\end{ex}

{\bf Proof of Corollary \ref{cor:main csp^2=csp}:}
This is a direct consequence of Theorem \ref{thm:main csp} and Theorem \ref{thm:main csp^2 cO^2} as indicated by Amdeberhan \cite{Amd12}.
By setting $t=0$ in Theorem \ref{thm:main csp} and $t_1=t_2=0$ in Theorem \ref{thm:main csp^2 cO^2},
we have
\begin{align*}
	\sum_{n=0}^\infty
	q^n
	\sum_{\lambda\vdash n}
	\prod_{u\in\lambda}
	\frac{c_{sp}(u)}{h(u)}
	=\prod_{n=1}^\infty \frac{1}{1+q^{4n-2}}
	=\sum_{n=0}^\infty
	\big((-1)^{1/2}q\big)^n
	\sum_{\lambda\vdash n}
	\prod_{u\in\lambda}
	\frac{c^2_{sp}(u)}{h^2(u)}.
\end{align*}
From the infinite product formula,
it is easy to see that,
in the above equation,
only terms labeled by $q^{n}, n\in2\mathbb{Z}_{\geq0}$ make non-trivial contribution.
Since $(-1)^{n/2}=(-1)^{\binom{n}{2}}$ when $n\in2\mathbb{Z}_{\geq0}$,
this corollary is obtained by taking the coefficient of $q^n$ on both sides of the above equation for an arbitrary fixed nonnegative integer $n$.
$\Box$

\ \\
\indent{\bf Proof of Theorem \ref{thm:main csp*c}:}
We will use the following Stanley's content formula \cite{S99} for Schur functions
\begin{align*}
	s_\lambda(1^n)
	=\prod_{u\in\lambda}
	\frac{n+c(u)}{h(u)}
\end{align*}
and the following Cauchy identity for symplectic Schur functions (see Theorem 15.1 in \cite{Sun86} and \cite{Sun90})
\begin{align}\label{eqn: cauchy for sps}
	\sum_{\lambda\in\mathcal{P}}
	q^{|\lambda|}
	sp_\lambda\big(\bm p(\bm y)\big) s_\lambda\big(\bm p'(\bm y')\big)
	=&\frac{\prod_{1\leq i<j}(1-q^2 y'_iy'_j)}
	{\prod_{i,j=1}^\infty(1-qy'_iy_j)}.
\end{align}
Actually,
Sundaram \cite{Sun86,Sun90} provided the above equation as an equality for formal functions in the variables $(x_1,x_1^{-1},...,x_n,x_n^{-1})$ for $l(\lambda)\leq n$.
We regard it as an identity in the ring $\mathbb{C}[\![q]\!][\bm p(\bm y), \bm p'(\bm y')]$.
For completeness,
in the following Lemma \ref{lem: cauchy for sps},
we provide a proof of the equation \eqref{eqn: cauchy for sps} using the method employed in deriving Theorem \ref{thm:sum spsp as p_k}.

By evaluating the above equation \eqref{eqn: cauchy for sps} at
\[y_i=y_j'=1\text{\ for\ }1\leq i \leq n_1, 1\leq j\leq n_2
\ \ \ \text{and}\ \ \ 
y_i=y_j'=0\text{\ for\ }i>n_1, j>n_2,\]
one obtains
\begin{align}
	\sum_{\lambda\in\mathcal{P}}
	q^{|\lambda|}
	\prod_{u\in\lambda}
	\frac{\big(n_1+c_{sp}(u)\big)\big(n_2+c(u)\big)}{h^2(u)}
	=\frac{(1-q^2)^{\binom{n_2}{2}}}{(1-q)^{n_1n_2}}.
\end{align}
Since both sides of the above equation could be viewed as elements in the ring $\mathbb{C}[\![q]\!][n_1,n_2]$,
and we already have proved the equality between them for all $(n_1,n_2)\in\mathbb{Z}_{>0}^2$,
the above equation still holds when replacing $n_1, n_2$ by any complex numbers $t_1, t_2$.
The proof is complete.
$\Box$

Below,
we give a proof of equation \eqref{eqn: cauchy for sps} using the method developed in this paper.
\begin{lem}\label{lem: cauchy for sps}
	We have
	\begin{align}
		\sum_{\lambda\in\mathcal{P}}
		q^{|\lambda|}
		sp_\lambda\big(\bm p(\bm y)\big) s_\lambda\big(\bm p'(\bm y')\big)
		=&\frac{\prod_{1\leq i<j}(1-q^2 y'_iy'_j)}
		{\prod_{i,j=1}^\infty(1-qy'_iy_j)}.
	\end{align}
\end{lem}
{\bf Proof:}
First,
from the Cauchy identity for Schur functions,
it is known
\begin{align*}
	\sum_{\lambda\in\mathcal{P}}
	s_{\lambda}(\bm p) \cdot \langle\lambda|
	=\langle0| \Gamma_+(\bm p).
\end{align*}
Then from equations \eqref{eqn:sum sp as skew s} and \eqref{eqn:sum slambda' G},
we have
\begin{align*}
	\begin{split}
		\sum_{\lambda\in\mathcal{P}}
		q^{|\lambda|}
		sp_\lambda\big(\bm p'(\bm y')\big) s_\lambda\big(\bm p(\bm y)\big)
		=\langle0|\Gamma_+\big(\bm p'(\bm y')\big)
		\cdot q^{L_0} \cdot \Gamma_-\big(\bm p(\bm y)\big)
		\cdot \exp(G)|0\rangle.
	\end{split}
\end{align*}
The right hand side of above equation can be computed using the commutation relations \eqref{eqn:gamma+- comm relation}, \eqref{eqn:L0 gamma+_} and \eqref{eqn: gamma+ G}.
It is equal to
\begin{align*}
	\frac{1}
	{\prod_{i,j=1}^\infty(1-qy'_iy_j)}
	\cdot \langle0|q^{L_0} \cdot \Gamma_-\big(\bm p(\bm y)\big)
	\cdot \exp(G)
	\cdot \prod_{i=1}^\infty
	\Gamma_-(qy_i')^{-1}\Gamma_+(qy_i')
	\cdot |0\rangle.
\end{align*}
Then from commutation relations \eqref{eqn:gamma+- comm relation}, \eqref{eqn:L0 gamma+_}, and Lemma \ref{lem: gamma+- G} again,
we can move the operators of the form $\Gamma_-(\cdot)$ to the left side.
The above equation is equal to
\begin{align*}
	\frac{\prod_{1\leq i<j}(1-q^2 y'_iy'_j)}
	{\prod_{i,j=1}^\infty(1-qy'_iy_j)}
	\cdot \langle0|\Gamma_-\big(\bm p(q\bm y)\big)
	\cdot \prod_{i=1}^\infty
	\Gamma_-(q^2y_i')^{-1}
	\cdot q^{L_0} \cdot  \exp(G)
	\cdot \prod_{i=1}^\infty \Gamma_+(qy_i')
	\cdot |0\rangle.
\end{align*}
This lemma then follows from the equation \eqref{eqn:gamma preserve vacuum} and $\langle0|\exp(G)|0\rangle=0$.
$\Box$

\begin{rmk}
	It is suggested by the anonymous referee that some other specialisations of identities derived in this paper may be interesting.
	For example, the following specialisation of symplectic Schur functions was considered in \cite{CS12},
	\begin{align}\label{eqn:hbar-eva}
		sp_{\lambda,2n}(s,s^{-1},s^3,s^{-3},\dots,s^{2n-1},s^{-2n+1})
		=\prod_{u\in\lambda}
		\frac{\langle2n+c_{sp}(u)\rangle_s}
		{\langle h(u)\rangle_s},
	\end{align}
	where $\langle i\rangle_s:=s^{i}-s^{-i}$.
	This kind of specialisation for Schur functions is known as (see, after some changes of notation, I.3.Example 3 in \cite{Mac})
	\begin{align}\label{eqn:hbar-eva s}
		s_{\lambda}(t^{-2m+1},t^{-2m+3},\dots,t^{2m-1},0,\dots)
		=\prod_{u\in\lambda}
		\frac{\langle2m+c(u)\rangle_t}
		{\langle h(u)\rangle_t},
	\end{align}
	where, similarly, $\langle i\rangle_t:=t^{i}-t^{-i}$.
	Then,
	applying the specialisations \eqref{eqn:hbar-eva} and \eqref{eqn:hbar-eva s} to the identity \eqref{eqn: cauchy for sps} gives
	\begin{align}\label{eqn:q-version}
		\begin{split}
			\sum_{\lambda\in\mathcal{P}}
			&q^{|\lambda|}
			\prod_{u\in\lambda}
			\frac{\langle2n+c_{sp}(u)\rangle_s
				\cdot \langle2m+c(u)\rangle_t}
			{\langle h(u)\rangle_s
				\cdot \langle h(u)\rangle_t}
			=\frac{\prod_{i=1}^{4m-3}(1-q^2t^{2i-4m+2})^{d_{i,m}}}
			{\prod_{i=1}^{2m}\prod_{j=1}^{2n}(1-qs^{2j-2n-1}t^{2i-2m-1})},
		\end{split}
	\end{align}
	where for $1\leq i\leq 4m-3$, define
	\begin{align*}
		d_{i,m}:=\begin{cases}
			\lfloor (i+1)/2 \rfloor, & 1\leq i\leq 2m-1,\\
			\lfloor (4m-1-i)/2 \rfloor, & 2m\leq i\leq 4m-3.
		\end{cases}
	\end{align*}
	For examples,
	When setting $n=m=1$,
	the both sides of equation \eqref{eqn:q-version} are equal to
	\begin{align*}
		1+\frac{(t^{2}+1) (s^{2}+1)}{s t} q&+\frac{(s^{4}+s^{2}+1) (t^{4}+t^{2}+1)}{s^{2} t^{2}} q^{2}\\
		&+\frac{(s^{2}+1) (s^{4}+1) (t^{2}+1) (t^{4}+1)}{s^{3} t^{3}} q^{3}
		+O(q^4).
	\end{align*}
	When setting $n=2, m=3$,
	the both sides of equation \eqref{eqn:q-version} are equal to
	\begin{align*}
		1+\frac{(s^{6}+s^{4}+s^{2}+1) (t^{10}+t^{8}+t^{6}+t^{4}+t^{2}+1)}{s^{3} t^{5}} q
		+O(q^2).
	\end{align*}
\end{rmk}

\section{Conflict of interest and data availability statement}
The author states that there is no conflict of interest, and
no datasets were generated or analysed during the current study.

\vspace{.2in}
{\em Acknowledgements}.
The author is grateful to the anonymous referee for many excellent suggestions.
The author would like to thank Professors Amdeberhan, Andrews, and Ballantine for their valuable discussions and their interest in this work.
The author would also like to express gratitude to Professors Xiaobo Liu, Jian Zhou, and Xiangyu Zhou for their encouragement.
The author is supported by the NSFC grants (No. 12288201, 12401079),
the China Postdoctoral Science Foundation (No. 2023M743717),
and the China National Postdoctoral Program for Innovative Talents (No. BX20240407).

\vspace{.2in}

\renewcommand{\refname}{Reference}
\bibliographystyle{plain}
\bibliography{reference}
\vspace{30pt} \noindent
\end{document}